\documentclass[journal,10pt,twocolumn]{IEEEtran}
\usepackage{cite}
\usepackage{amsmath}
\usepackage{amsmath}
\usepackage{amsfonts}   
\usepackage{amssymb}    
\usepackage{float}
\usepackage{subfigure}
\usepackage{float}
\usepackage{amsmath}
\usepackage{color}
\usepackage{graphicx}
\usepackage{epsfig}
\usepackage{algorithm}
\usepackage{enumerate}
\usepackage[noend]{algorithmic}
\usepackage[singlespacing]{setspace}
\newcommand{\figref}[1]{Fig.~\ref{#1}}
\begin{document}
\title{A Framework for Uplink Intercell Interference Modeling with Channel-Based Scheduling\footnote{\thanks{$^\dag$ Part of this work has been presented at the IEEE International Symposium on Wireless Communication Systems (ISWCS 2011), Aachen,
Germany, November 2011.
 Hina Tabassum, Ferkan Yilmaz and Mohamed Slim Alouini are with the Computer, Electrical, Mathematical Sciences and Engineering (CEMSE) Division, KAUST, Thuwal, Makkah Province, Saudi Arabia. Email: \{hina.tabassum, ferkan.yilmaz, slim.alouini\}@kaust.edu.sa. Zaher Dawy is with the Department of Electrical and Computer Engineering, American University of Beirut (AUB), Beirut, Lebanon. Email: \{zaher.dawy\}@aub.edu.lb.
}
}}
{
\author{\IEEEauthorblockN{Hina Tabassum, \emph{Student Member, IEEE}, Ferkan Yilmaz, \emph{Member, IEEE}, Zaher Dawy, \emph{Senior Member, IEEE}, and Mohamed-Slim
Alouini, \emph{Fellow, IEEE}} 
}
}
\maketitle
\begin{abstract}
This paper presents a  novel framework for modeling the uplink intercell interference (ICI) in a multiuser cellular network. The proposed framework assists in quantifying the impact of various fading channel models and state-of-the-art scheduling schemes on the uplink ICI. Firstly, we derive a semi-analytical expression for the distribution of the location of the scheduled user in a given cell considering a wide range of scheduling schemes. Based on this, we derive the distribution and  moment generating function (MGF) of the uplink ICI considering a single interfering cell. Consequently, we determine the MGF of the cumulative ICI observed from all interfering cells and derive explicit MGF expressions for three typical fading models. Finally, we utilize the obtained expressions to evaluate important network performance metrics such as the outage probability, ergodic capacity, and average fairness numerically. Monte-Carlo simulation results are provided to demonstrate the efficacy of the  derived analytical expressions.
\end{abstract}
\IEEEpeerreviewmaketitle

\section{Introduction}

Explosive growth in the demand of high quality wireless data services compel the  network  designers to utilize spectrum more aggressively which on one side enhances the spectrum efficiency, whereas on the other side it enhances the
intercell interference (ICI) which  is  an alarming bottleneck in the telecommunication growth paradigm. 
The allocation of the same frequency bands across neighboring cells produces indeterministic ICI which is highly dependent on the statistics of the channel characteristics and on the dynamics of the  scheduling decisions. 
In this context,  it is of immense importance for the system designers to accurately characterize and investigate the behavior of the ICI  which  helps in gaining more theoretical insights,  quantifying various network performance metrics and developing efficient resource allocation and interference mitigation schemes.

Orthogonal frequency division multiple access (OFDMA) has been recently adopted as the multiple access scheme for the state-of-the-art LTE and WiMAX cellular technologies. In OFDMA, a wide-band frequency-selective fading channel is decomposed into a set of orthogonal narrow-band  subcarriers. The orthogonality among the subcarriers per cell makes the intra-cell interference almost negligible. However,  with universal frequency reuse among cells (i.e., all cells use the same set of subcarriers), the ICI at each subcarrier may cause severe degradation in the network performance. In OFDMA networks, the subcarriers are allocated adaptively among users per cell based on a predefined scheduling scheme. Moreover, each subcarrier is allocated to only one user per cell assuming BSs are equipped with single antenna
and, thus, the number of interfering users on each subcarrier is rather limited. Therefore, the cumulative ICI on a given subcarrier may not be modeled accurately as a Gaussian random variable (RV) by simply invoking the central limit theorem.

Several recent studies considered the modeling of ICI in the downlink where the location of interferers is typically deterministic.
A semi-analytical distribution for the signal-to-interference-noise ratio (SINR) has been derived in \cite{eurasip} under path loss and log normal shadowing for randomly located femtocell networks. In \cite{plass}, the applicability of the Gaussian and binomial distributions for modeling the downlink ICI is investigated. In \cite{yang},  the optimal threshold is derived for fractional frequency reuse (FFR) systems  assuming ICI as Gaussian RV. In \cite{TWC1}, the authors derived the distribution of the ICI under log-normal shadowing and Rayleigh fading. The distribution of ICI is shown to highly deviate from the Gaussian distribution in OFDMA networks.

In comparison to downlink, the nature of uplink ICI is different in various aspects that include the following: 
(i) Due to the implicit symmetry and fixed locations of the BSs in the typical grid-based downlink network models, the number of significantly contributing interferers typically remains the same irrespective of the position of the mobile receiver. Also, it has been shown in \cite{plass} that the strongest interference is generated by  two closest interfering BSs  irrespective of the mobile receiver location. However, the number of significantly contributing interferers in the uplink cannot be quantified at a given instant due to the highly varying locations of the interfering mobile transmitters;
(ii) Conditioned on the location of the desired mobile receiver within a cell, the exact  distance of the interfering BSs can be calculated in the typical grid-based downlink network models. However, knowing  the location of the BS receiver in the uplink  does not help in determining the exact location of the interfering mobile users;
(iii) In the uplink, cell edge and cell center mobile users are subject to the same amount of interference on a given subcarrier, which is the interference received at the BS.  Whereas the same is not true for the downlink  in which cell edge  users  experience higher interference coming from the nearby  BSs \cite{TWC,viering}.

Some worth mentioning  research works for the uplink appear in \cite{TWC, IEEEletter,rup,ICC}.
In \cite{TWC}, the authors developed an analytical model for subcarrier collisions as a function of the cell load and frequency reuse pattern. They derived an expression for the SINR in the uplink and downlink, ignoring the effect of shadowing and fading. In \cite{IEEEletter}, the authors developed an analytical expression for the subcarrier collision probability considering non-coordinated schedulers. In \cite{rup}, the authors modeled uplink ICI in an OFDMA network as a function of the reuse partitioning radius and traffic load assuming arbitrary scheduling. In \cite{ICC}, the authors presented a semi-analytical method to approximate the distribution of the uplink ICI through numerical simulations without considering the impact of scheduling schemes.

In this paper, we propose a novel theoretical framework to derive the statistics of  the uplink ICI  on a given subcarrier as a function of both the channel statistics (i.e., path loss, shadowing and fading) and scheduling decisions. The proposed framework can be also extended to typical downlink scenarios  as explained in \cite{hina}. The framework is generic in the sense that the derivations hold for generalized fading channels and various scheduling algorithms.   We start by deriving the distribution of the location of the scheduled user in a given cell. We then derive the distribution and moment generating function (MGF) of the ICI considering a single interfering cell. Next, we derive the MGF expression for the cumulative ICI experienced from all interfering cells over generalized fading channels, and present explicit expressions for three practical fading models. Finally, we demonstrate the importance of the derived expressions by utilizing them to evaluate important network performance metrics \textcolor{black}{such as outage probability and ergodic capacity.} 

The remainder of this paper is organized as follows. Section~II presents the system model and the main steps of the proposed framework. In Section~III, the distribution of the scheduled user location is derived for different scheduling algorithms. In Section~IV, the distribution of the uplink ICI from one neighboring cell is derived. The MGF of the cumulative ICI from all interfering cells is determined in Section~V and utilized in Section~VI to evaluate three network performance metrics. Finally, numerical and simulation results are presented and analyzed in Section~VII, and conclusions are drawn in Section~VIII.

\noindent\textbf{Notation}:  $\mathrm{Exp}(\lambda)$ represents an exponential distribution with parameter $\lambda$, $\mathrm{Gamma}(m_s,m_c)$ represents a Gamma distribution with shape parameter $m_s$ and scale parameter $m_c$. $\mathcal{K}_G (m_c,m_s,\Omega)$ represents the Generalized-$\mathcal{K}$ distribution with  fading parameter $m_c$, shadowing parameter $m_s$ and average power $\Omega$. $\Gamma(.)$ represents the Gamma function. $\textcolor{black}{P(A)}$ denotes the probability of event A. $f(.)$ and $F(.)$ denotes the probability distribution function (PDF) and cumulative distribution function (CDF), respectively. $[a,b]$ denotes a discrete set of elements which ranges from $a$ to $b$. Finally, $\mathbb{E}[.]$ denotes the expectation operator.

\begin{figure}[t]
  \centering
  \includegraphics[scale=0.7]{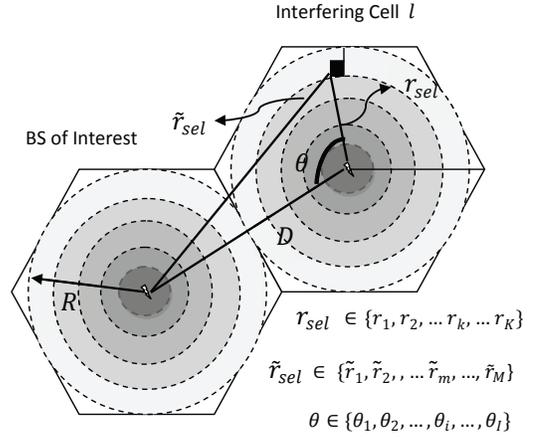}
\caption{Geometrical illustration of dividing the cellular network into multiple rings of non-uniform width $\Delta_k$.}
\label{Cap0}
\vspace{-7mm}
\end{figure}
\section{System Model and Proposed Framework}
\subsection{Description of the System Model}
We consider a given cell surrounded by $L$ interfering neighboring cells. For analytical convenience, the cells are assumed to be circular with radius $R$. Each cell $l$ is assumed to have $U$ uniformly distributed users. 
The frequency reuse factor is assumed to be unity with each subcarrier reused in all cells.
The bandwidth of a subcarrier is assumed to be less than the channel coherence bandwidth, thus, each subcarrier experiences flat fading. Time is divided into time slots of length smaller than the channel coherence time and, thus, the channel variation within a given time slot is negligible.

Generally, the scheduling strategies can be broadly categorized into two classes; (i) rate maximization (i.e., rate adaptation) while transmitting with  constant/maximum power; (ii) power minimization (i.e., power adaptation) while achieving a fixed data rate. In this work,  \emph{we focus on rate adaptive schemes where users transmit with their maximum power in order to maximize their rate depending on the existing channel and interference conditions.}
Therefore, for the scope of this paper, we assume that all users transmit with their maximum power $P_{\mathrm{max}}$ on a given subcarrier with rate adaptation depending on their channel qualities. At this point, it is however important to emphasize that this is not a limitation and the approach can be extended for  uplink power control mechanisms as discussed in \cite{hina}.  The instantaneous signal to noise ratio (SNR) $\gamma$ of each user  can then be written as follows:
\begin{equation}
\label{1}
\gamma=P_{\mathrm{max}} C \frac{{r}^{-\beta} \psi {\eta} }{\sigma^2}= \bar{K} {{r}^{-\beta} \zeta},
\end{equation}
where $\bar{K}=\frac{P_{\mathrm{max}} C}{\sigma^2}$,
$C$ is the path loss constant, $r$ is the user distance from its serving BS, $\psi$ and $\eta$ denotes the shadowing and small scale  fading coefficient between user and BS on a given subcarrier, respectively, $\beta$ is the path loss exponent, $\sigma$ denotes the thermal noise at the receiver and $\zeta$ is the composite fading. Note that all users are assumed to be associated with their closest BS \cite{TWC,novlan}, therefore $r \leq R$. 

Each cell is divided into $K$ concentric circular regions. Since path loss decays exponentially from cell center to cell edge, therefore, we consider discretization of cellular region in such a way that the path loss decay within each circular region remains constant or uniform. The main motivation for dividing the cell into a discrete set of circular regions relies on the fact that the channel statistics of the users located in a given circular region become relatively similar especially for large values of $K$.  More explicitly, the characterization of the circular regions can be demonstrated as follows:
\begin{equation}
\mathrm{log}_{10}r_k={\frac{\kappa+ 10 \beta \mathrm{log}_{10} r_{k-1}}{10 \beta}},\:\:\: r\leq R,
\end{equation}
where $\kappa$ is the path loss decay within each circular region [dB]. Due to the exponential nature of the path loss, it varies rapidly near the cell center than at the cell edge, therefore, (i) each of the $k^{\mathrm{th}}$ circular region bounded by two adjacent rings, i.e., $r_k$ and $r_{k-1}$ possess non-uniform width $\Delta_k=r_k-r_{k-1}$; (ii) the number of circular regions are high in the cell center than at the cell edge; (iii)  the average number of users located within $k^{\mathrm{th}}$ circular region bounded by ring $r_k$ and $r_{k-1}$ are considered to be located at $r_k$. Note that, this is an approximation which is required for deriving the analytically tractable model of ICI and in any case it is not required for the Monte-Carlo simulations. The average number of users in each ring $k$ (for analysis) can then be given as:
\begin{equation}
u_k=\frac{U(r_{k}^2-r_{k-1}^2)}{R^2}, \:\:\:\:\: k=1,2,\cdots,K,
\end{equation}
It is important to note that $u_k$ can be a fraction; therefore, we round off the fractional part of users.
\subsection{Main Steps of the Proposed Framework}
In order to characterize the statistics of the  uplink ICI for generalized fading channels and  various scheduling schemes, the proposed framework mandates the following steps:
\begin{enumerate}[i)]
\item Derive the distribution $f_{r_{\mathrm{sel}}}(r)$ of the distance of the allocated user $r_{\mathrm{sel}}$ in a given cell from its serving BS based on the deployed scheduling scheme. 
\item Derive the distribution $f_{\tilde{r}_{\mathrm{sel}}}(\tilde{r})$ of the distance  between the allocated user in a neighboring interfering cell and the BS of the cell of interest $\tilde{r}_{\mathrm{sel}}$.
\item Derive the distribution $f_{X_l}(x)$ of the interference from the allocated user in neighboring cell $l$ to the BS of interest, i.e., $X_l=\bar{K}\tilde{r}_{\mathrm{sel}}^{-\beta}\chi$,
where $\chi$ denotes the combined shadowing and fading, i.e., composite fading of the interference. Finally, derive the MGF of the cumulative ICI, i.e., ${Y}=\sum_{l=1}^{L}X_l$ from the scheduled users on a given subcarrier in all $L$ interfering cells.
\end{enumerate}

\section{Distribution of the Scheduled User Location}
Considering the dependence of the uplink ICI on the location of the scheduled users in the neighbor interfering cells which  in turn depends on the deployed scheduling schemes,  we derive in this section the distribution of the distance between the scheduled user and its serving BS in a given cell (i.e., the probability mass function (PMF) of $r_{\mathrm{sel}}$) considering the following five scheduling algorithms: greedy scheduling, proportional fair scheduling, round robin scheduling, location based round robin scheduling, and greedy round robin scheduling.

\subsection{Greedy Scheduling Scheme}
\noindent  Greedy scheduling is an opportunistic scheme that aims at maximizing the network throughput by taking full advantage of multiuser diversity. However, it suffers from low fairness among users which makes it less attractive for network operators.  The procedure for determining the PMF of $r_{\mathrm{sel}}$ considering greedy scheduling is divided into two steps:\\
\noindent\textbf{Step~1~(\rm\textit{Selecting the user with the highest SNR in ring $k$}):}
Since the path loss decay within each circular region is considered to be uniform, we approximate the distance of all  users located within $k^{\mathrm{th}}$ circular region by ring $r_k$ for analytical tractability as we already mentioned in Section II. In this step, we select a user with maximum SNR in each ring $k$ which  possess $u_k$ users. Thus, selecting a user in a ring $k$ is equivalent to selecting the user with maximum channel gain among all the users in ring $k$, i.e.,
$
\zeta_k=\text{max}\{\zeta_1,\zeta_2,\cdots,\zeta_i,\cdots, \zeta_{u_k}\}
$,
where $\zeta_i$ is the composite fading channel gain  between user $i$ and its BS on a given subcarrier. 
The CDF and PDF of the maximum channel gain $\zeta_k$ can be written as follows:
\begin{equation}
\label{CDF1}
\begin{split}
F_{\zeta_k}(\zeta_k)=\prod_{i=1}^{u_k}F_{\zeta_i}(\zeta_k)\stackrel{\mathrm{i.i.d}}{=}\left(F_{\zeta}(\zeta_k)\right)^{u_k},
\end{split}
\end{equation}
\begin{equation}
\label{PDF1}
f_{\zeta_k}(\zeta_k)=\sum_{j=1}^{u_k} f_{\zeta_j}(\zeta_k)\prod_{i=1,i\neq j}^{u_k}F_{\zeta_i}(\zeta_k)\stackrel{\mathrm{i.i.d}}{=}{u_k}f_{\zeta}(\zeta_k)\left(F_{\zeta}(\zeta_k)\right)^{u_k-1},
\end{equation}
To consider path loss, we now perform a transformation of RVs using \eqref{1}, $\gamma_{k}=\bar{K}r_{k}^{-\beta}\zeta_k$, where, $\gamma_{k}$ is the selected user SNR in each ring $k$. The CDF and PDF of $\gamma_{k}$ can then be written as follows:
\begin{equation}
\label{CDFstep1}
F_{\gamma_{k}}(\gamma_{k})= \prod_{i=1}^{u_k}F_{\zeta_i}(\bar{K}^{-1}\gamma_{k}r_{k}^{\beta} ) \stackrel{\mathrm{i.i.d}}{=}\left(F_{\zeta}(\gamma_{k}r_{k}^{\beta} \bar{K}^{-1})\right)^{u_k},
\end{equation}
\begin{equation}
\label{PDFstep1}
\begin{split}
f_{\gamma_{k}}(\gamma_{k})=\frac{1}{\bar{K} r_k^{-\beta}} \sum_{j=1}^{u_k} f_{\zeta_j}(\gamma_{k}r_{k}^{\beta}\bar{K}^{-1} )\prod_{i=1,i\neq j}^{u_k}F_{\zeta_i}(\gamma_{k}r_{k}^{\beta} \bar{K}^{-1}),\\\stackrel{\mathrm{i.i.d}}{=}\frac{u_k}{\bar{K} r_{k}^{-\beta}}f_{\zeta}(\gamma_{k} r_{k}^{\beta}\bar{K}^{-1})\left(F_{\zeta}(\gamma_{k} r_{k}^{\beta}\bar{K}^{-1})\right)^{u_k-1}.
\end{split}
\end{equation}

\noindent\textbf{Step~2~(\rm\textit{Selecting the user with maximum SNR among $K$ rings}):}
In this step, we  compute the probability of selecting the $k^\mathrm{th}$ ring among all other rings. It is important to note that this is equivalent to selecting the ring $k$ which possesses the user with the highest SNR among all rings.
Conditioning on $\gamma_{k}$, the PDF of $r_{\mathrm{sel}}$  can be written explicitly as follows:
\begin{equation}
\label{Final0}
P(r_{\mathrm{sel}}=r_{k}|\gamma_{k})=\prod_{i=1,i\neq k}^{K}p(\gamma_{i}\leq \gamma_{k})=\prod_{i=1,i\neq k}^{K} F_{\gamma_{i}}(\gamma_{k}),
\end{equation}
By averaging over the distribution of $\gamma_{k}$, the final expression for the PMF of $r_{\mathrm{sel}}$ is
\begin{equation}
\label{Final}
P(r_{\mathrm{sel}}=r_k)=\int_0^\infty\left({\prod_{i=1,i\neq k}^{K}F_{\gamma_{i}}( \gamma_{k})}\right)f_{\gamma_{k}}(\gamma_{k}) d\gamma_{k},
\end{equation}
Using \eqref{CDFstep1}, \eqref{Final} can be written for i.i.d. case as follows:
\begin{equation}
\label{Final2}
\begin{split}
P(r_{\mathrm{sel}}=r_{k})=\int_0^\infty\prod_{i=1,i\neq k}^{K}
\left(F_{\zeta}(\gamma_{k}r_{i}^{\beta}\bar{K}^{-1})\right)^{u_i}
 \times \\
\frac{u_k f_{\zeta}(\gamma_{k}r_{k}^{\beta}\bar{K}^{-1})}{\bar{K}r_{k}^{-\beta}}\left(F_{\zeta}(\gamma_{k} r_{k}^{\beta}\bar{K}^{-1} )\right)^{u_k-1}.
\end{split}
\end{equation}
where $r_{\mathrm{sel}} \in [0,R]$.
The results in \eqref{Final2} are generalized for any shadowing and fading statistics. Even though \eqref{Final2} is not a closed form expression, the integration can be solved accurately and efficiently using standard mathematical software packages such as \texttt{MAPLE} and \texttt{MATHEMATICA}.

\subsection{Proportional Fair Scheduling Scheme}
The proportional fair scheduling scheme  allocates the subcarrier to the user with the largest normalized SNR (${\gamma}/{\bar{\gamma}}$)  \cite{prop1}, where $\gamma$ and $\bar{\gamma}$ denote the instantaneous SNR and the short term achieved average SNR of a given user, respectively. In other words, the selection criterion is based on selecting a user who has maximum instantaneous SNR relative to its own average SNR.
The distribution of $r_{\mathrm{sel}}$ can be derived as:\\
\noindent\textbf{Step~1~(\rm\textit{Selecting the user with maximum normalized SNR in ring $k$}):}
In this step, the performance of proportional fair scheduling scheme is independent of the path loss factor if users are moving relatively slowly, i.e., their path loss remains nearly the same on a short term basis. 
In this case, the problem of selecting the maximum normalized SNR in a ring $k$ can be written as:
\begin{equation}
\label{Norm1}
{\zeta_k}
=\text{max}\left\{\frac{\zeta_1}{\bar{\zeta}_1},\frac{\zeta_2}{\bar{\zeta}_2},\cdots,\frac{\zeta_i}{\bar{\zeta}_i},\cdots,\frac{\zeta_{u_k}}{\bar{\zeta}_{u_k}}\right\},\:\:\:
\end{equation}
where $\bar{\zeta}_i=\int_0^\infty \zeta_i f_{\zeta_i}(\zeta_i) d\zeta_i$ is the average of the composite fading channel and $\zeta_k$ is the maximum normalized composite fading channel gain  in ring $k$.
For i.i.d. average composite fading gains  of the users located in ring $k$, i.e.,  $\bar{\zeta}=\bar{\zeta}_1=\bar{\zeta}_2=\cdots \bar{\zeta}_{u_k}$, the problem of selecting the user with  maximum normalized channel gain reduces to selecting the user with the maximum channel gain, i.e.,
$
\zeta_k=\text{max}\left\{{\zeta_1},{\zeta_2},\cdots,{\zeta_i},\cdots,{\zeta_{u_k}}\right\}
$.
Thus, for any ring $k$, the CDF and PDF of the selected SNR  $\gamma_k=\bar{K}r_k^{-\beta}\zeta_k$ can be written as:
\begin{equation}
F_{\gamma_k}(\gamma_k)=
\prod_{i=1}^{u_k}F_{\zeta_i}(\bar{\zeta_i}\gamma_k r_k^{\beta}\bar{K}^{-1})\stackrel{\mathrm{i.i.d}}{=}(F_{\zeta_i}(\bar{\zeta}\gamma_k r_k^{\beta}\bar{K}^{-1}))^{u_k},
\end{equation}
\begin{equation}
\begin{split}
f_{\gamma_k}(\gamma_k)=\frac{r_k^{\beta}}{\bar{K}} \sum_{j=1}^{u_k} f_{\zeta_j}(\bar{\zeta_j}\gamma_k r_k^{\beta}\bar{K}^{-1})\prod_{i=1,i\neq j}^{u_k}F_{\zeta_i}(\bar{\zeta_i}\gamma_k r_k^{\beta}\bar{K}^{-1}),\\
\stackrel{\mathrm{i.i.d}}{=}\frac{u_k r_k^{\beta}}{\bar{K}} (F_{\zeta_i}(\bar{\zeta}\gamma_k r_k^{\beta}\bar{K}^{-1}))^{u_k-1}f_{\zeta}(\bar{\zeta}\gamma_k r_k^{\beta}\bar{K}^{-1}).
\end{split}
\end{equation}
The short term average SNR of the selected user in ring $k$, i.e., $\bar{\gamma}_k$ can then be computed as $
\bar{\gamma}_k=\int_0^\infty \gamma_k f_{\gamma_{k}}(\gamma_k)d\gamma_k
$. Finally, the normalized selected SNR in each ring $k$ can be defined as $\xi_k=\frac{\gamma_k}{\bar{\gamma}_k}$ and performing a transformation of RVs, the CDF and PDF of $\xi_k$ can be given as 
$
F_{\xi_k}(\xi_k)=F_{\gamma_k}(\bar{\gamma}_k\xi_k)
$; and
$
f_{\xi_{k}}(\xi_k)={\bar{\gamma}_k}f_{\gamma_k}(\bar{\gamma}_{k} \xi_k)
$, respectively.
\\
\noindent\textbf{Step~2~(\rm\textit{Selecting the ring $k$ with maximum normalized SNR from the $K$ rings}):} Once we characterize the PDF and CDF of $\xi_k$, the probability of selecting any ring $k$ can be written using \eqref{Final0} as follows:
\begin{equation}
\label{Finalnorm}
P(r_{\mathrm{sel}}=r_k)=\int_0^\infty\left({\prod_{i=1,i\neq k}^{K}F_{\xi_{i}}( \xi_{k})}\right)f_{\xi_{k}}(\xi_{k}) d\xi_{k}.
\end{equation}

\subsection{Round Robin Scheduling Scheme}
Round robin scheduling is a non-opportunistic scheme where a user is selected randomly within a time slot. As each user has equal probability of allocation, round robin is a \emph{strictly fair} scheduling scheme. The round robin scheme provides maximum fairness among users and may serve as a lower bound in terms of network throughput which is useful in calibrating the performance of other scheduling schemes, however, the resulting network throughput is significantly low which makes it less attractive for practical implementations.  The PMF of the scheduled user location can then be given as
\begin{equation}
P(r_{\mathrm{sel}}=r_{k}) = \frac{u_k}{U}.
\end{equation}

\subsection{Location Based Round Robin Scheduling Scheme}
Location based round robin is another non-opportunistic scheduling scheme which do not require any channel state information, however, it requires the location information of the users. Even though location based scheduling is not common in practice,
the location of each mobile user can be determined  at the BS using global positioning system (GPS) or estimate based on a power measurement of pilot signals from the surrounding beacons, e.g., using triangulation based techniques. In this regard, there are variety of techniques available in the literature which demonstrate  how the location of users can be evaluated at the BS (see\cite{loc1,loc2} and the references therein). 
Moreover, the users located in different circular regions can also be classified based on the long term average SNRs, i.e., by computing SNR thresholds for different distances (rings) \cite{yang} which is a common technique in FFR systems to distinguish between cell-edge and cell center users.

In this scheme, we consider $W$ time slots during which the distance of the users from their serving BS will remain approximately the same. For simplicity, the number of time slots $W$ is set equal to $K$. 
At a given time slot $T_w$, we select any arbitrary user from a specific ring (analysis) and circular region (simulations) starting from the cell center. We continue to allocate the users by accessing the circular rings sequentially from cell center to cell edge. At this point, it is important to emphasize that all cells are considered to be time synchronous in allocating the users from particular rings, i.e., at a given time slot $T_w$ all cells are selecting the $w^{th}$ ring. Thus, the    PMF of $r_{\mathrm{sel}}$ for a given time slot $w$ denoted by $P(r_{\mathrm{sel}}=r_k^{T_{w}})$ can be given as:
\begin{equation}
P(r_{\mathrm{sel}}=r_k^{T_w})=
\begin{cases}
1, & \text{if }k = w\\
0, &  \text{else}
\end{cases}.
\end{equation}

\subsection{Greedy Round Robin  Scheduling Scheme}
Greedy round robin is an opportunistic scheduling scheme which captures the multiuser diversity  while maintaining some fairness among users. In this scheme,
we consider $W=K$ time slots during which the distance of the users from  their serving BSs remain nearly the same, however, the small scale fading gain on the considered subcarrier may vary from one time slot to the other. We select the user with maximum SNR  in each time slot $T_w$, however once a user is selected from a ring, all users located in that ring will not be scheduled for transmission for the next $K-1$ time slots.
Note that all BSs are considered to be time synchronized in terms of scheduling.
Clearly, the  probability of allocating a ring $k$ at $T_1$ can simply be given by \eqref{Final2}. However, the probability of selecting a ring $k$ at $T_2$ is a dependent event and can be derived using Bayes theorem as follows:
\begin{equation}
P(r_{\mathrm{sel}}=r_k^{T_{2}})=\sum_{\substack{
j=1\\
j \neq k}}^K P\left(\frac{{r_{\mathrm{sel}}=r_k}^{T_2}}{{r_{\mathrm{sel}}=r_j}^{T_1}}\right)P(r_{\mathrm{sel}}=r_j^{T_1}),
\end{equation}
where,
\begin{equation}
\begin{split}
P\left(\frac{{r_{\mathrm{sel}}=r_k}^{T_2}}{{r_{\mathrm{sel}}=r_j}^{T_1}}\right)=
\int_0^\infty\prod_{\substack{i\neq k\\ i\neq j}}^{K}\left(F_{\zeta}(\gamma_{k}r_{i}^{\beta}\bar{K}^{-1})\right)^{u_i}\times\\ \frac{u_k f_{\zeta}(\gamma_{k}r_{k}^{\beta}\bar{K}^{-1})}{\bar{K} r_{k}^{-\beta}}\left(F_{\zeta}(\gamma_{k} r_{k}^{\beta}\bar{K}^{-1} )\right)^{u_k-1}d\gamma_{k}.
\end{split}
\end{equation}
Since the probability of allocating any ring $k$ within time slot $T_w$ depends on all previous states, therefore, the principle of Markov chain transition probabilities is not directly applicable. For more clarity,  the  probability of selecting a ring $k$ at $T_3$ is given as follows:
\begin{equation}
\begin{split}
P(r_{\mathrm{sel}}=r_k^{T_{3}})=\sum_{m \neq k}^{K}\sum_{\substack{j \neq m\\ j\neq k}}^K P\left(\frac{{r_{\mathrm{sel}}=r_k}^{T_3}}{{r_{\mathrm{sel}}=r_j}^{T_2}\cap {r_{\mathrm{sel}}=r_m}^{T_1}}\right)\times \\ P\left(\frac{{r_{\mathrm{sel}}=r_j}^{T_2}}{{r_{\mathrm{sel}}=r_m}^{T_1}}\right)P(r_{\mathrm{sel}}=r_m^{T_1}),
\end{split}
\end{equation}
\textcolor{black}{where, $P\left(\frac{{r_{\mathrm{sel}}=r_k}^{T_3}}{{r_{\mathrm{sel}}=r_j}^{T_2}\cap {r_{\mathrm{sel}}=r_m}^{T_1}}\right)$ is given at the top of next page in \eqref{top}.}
\begin{figure*}[t]
\begin{equation}\label{top}
\begin{split}
P\left(\frac{{r_{\mathrm{sel}}=r_k}^{T_3}}{{r_{\mathrm{sel}}=r_j}^{T_2}\cap {r_{\mathrm{sel}}=r_m}^{T_1}}\right)=
\int_0^\infty\prod_{\substack{i\neq k \\ i\neq j, i\neq m}}^{K}\left(F_{\zeta}(\gamma_{k}r_{i}^{\beta}\bar{K}^{-1})\right)^{u_i}
\frac{u_k f_{\zeta}(\gamma_{k}r_{k}^{\beta}\bar{K}^{-1})\left(F_{\zeta}(\gamma_{k} r_{k}^{\beta}\bar{K}^{-1} )\right)^{u_k-1}}{\bar{K}r_{k}^{-\beta}}d\gamma_{k},
\end{split}
\end{equation}
\hrule
\end{figure*}
In general, the probability of selecting any ring $k$ at a time slot $T_w$, i.e.,  $P(r_{\mathrm{sel}}=r_k^{T_{w}})$ can be written as:
\begin{multline}
\sum_{n \neq k}^{K}\sum_{\substack{s\neq n, \\s\neq k}}^K \cdots \sum_{\substack{j\neq s,n,..\\j\neq k}}
P(r_{\mathrm{sel}}=r_n^{T_1})
P\left(\frac{{r_{\mathrm{sel}}=r_s}^{T_2}}{{r_{\mathrm{sel}}=r_n}^{T_1}}\right)
\times \cdots \\\times
P\left(\frac{{r_{\mathrm{sel}}=r_k}^{T_w}}{{r_{\mathrm{sel}}=r_j}^{T_{w-1}}\cdots \cap {r_{\mathrm{sel}}=r_s}^{T_2}\cap {r_{\mathrm{sel}}=r_n}^{T_1} }\right). 
\end{multline}
{\it Computational Efficiency:} The time complexity of the greedy round robin scheme is  heavily based on the computational time of the  \texttt{NIntegrate} operation in \texttt{MATHEMATICA}. The computation time of one \texttt{NIntegrate} operation is denoted by $\tau$ where $\tau =$ 0.95 sec (on Intel(R) Xeon(R) CPU-X5550 @2.67 GHz with 24 GB RAM, 64 bit Operating System); this is equivalent to the computational complexity of (i) \emph{greedy} scheme and (ii) the \emph{first time slot of the greedy round robin} scheme.  Monte-Carlo simulation time required for 100,000 trials in \emph{greedy} scheme requires around 150.67~sec which is  150 times more than the computational complexity of \texttt{NIntegrate} operation. This fact demonstrates the computational efficiency of \emph{greedy} scheme in comparison to Monte-Carlo simulations. 
However, in the second  time slot, greedy round robin scheme requires $K-1$ integrations whereas, for the third time slot $K-2$ integrations are required. Therefore, the computation time at any time slot $w \leq W$ can be given as $\tau+\sum_{i=2}^{w} \tau(K-i+1);\:\:\:w\leq W$, where $W$ denotes the total number of time slots. Considering $w = W = $15 and the number of rings $K=W$  for greedy round robin,  the analytical time complexity  is around 113~sec which is still lower than the Monte-Carlo simulation time required for the \emph{greedy} scheme. Therefore, it can be concluded  that the numerical implementation of greedy round robin would be beneficial for performance assessment in cases with reduced number of rings, e.g, when either the number of slots, the path-loss exponent, or the cell radius has a small value.

\subsection{Evaluating the Joint PMF of $r_{\mathrm{sel}}$ and $\theta$}
Note that, $P(r_{\mathrm{sel}}=r_{k})$ derived for all of the above scheduling schemes is the marginal PMF of $P(r_{\mathrm{sel}}=r_k,\mathcal{\theta}=\theta_i)$ where $\mathcal{\theta}$ denotes the angle of the allocated user with respect to the serving BS and it is uniformly distributed from $0$ to $2\pi$ (see Fig. 1). Although the PDF of $\mathcal{\theta}$ is continuous we can discretize it for analytical consistency and complexity reduction.
Consider discretizing the range of RV $\theta$ in ${\mathcal{I}}$ uniform angular intervals of desired accuracy. Thus $P(\mathcal{\theta}=\theta_i)=\frac{1}{{\mathcal{I}}}$ where $\theta_i$ denotes any discrete value that the RV $\mathcal{\theta}$ can take. Since $r_{\mathrm{sel}}$ and $\mathcal{\theta}$ are independent, their joint PMF can be written as:
\begin{equation}
\label{IntPMF}
\begin{split}
P(r_{\mathrm{sel}}=r_k, \mathcal{\theta}=\theta_i)=\frac{P(r_{\mathrm{sel}}=r_k)}{\mathcal{I}}.
\end{split}
\end{equation}

\section{Distribution of Intercell Interference from One Cell}
The derivation for the distribution of the ICI from an interfering cell $l$, i.e., $f_{X_l}(x)$, depends on the distribution of the distance between the allocated user in the interfering cell $l$ and the BS of interest, i.e., $f_{\tilde{r}_{sel_l}}(\tilde{r})$.  As mentioned earlier, each interfering cell is assumed to have identical conditions in a given time slot. Therefore, $f_{\tilde{r}_{\mathrm{sel}}}(\tilde{r})$ applies to all interfering cells and, thus, we will drop the subscript $l$ in the sequel to simplify notation. Using the cosine law (see \figref{Cap0}), we can write:
\begin{equation}
\label{cosine law}
\begin{split}
\tilde{r}_{{\mathrm{sel}}}^2=r_{{\mathrm{sel}}}^2+D^2-2{r}_{{\mathrm{sel}}}D \;\mathrm{cos}{\theta}.
\end{split}
\end{equation}
$\tilde{r}_{{\mathrm{sel}}}$ is the distance of the allocated user in the interfering cell $l$ from the BS of interest, $r_{{\mathrm{sel}}}$ is the distance of the allocated user from its serving BS, i.e., (BS $l$), $\theta \in [0, 2\pi]$ and $D = 2R$ since we consider universal frequency reuse with one tier of interfering cells. The approach can be extended to any number of tiers in a straightforward manner.
In order to determine the PMF of $\tilde{r}_{\mathrm{sel}}$ where $ \tilde{r}_{\mathrm{sel}} \in [D-R, D+R] $, first of all we compute $\tilde{r}_{i,k}$ for given $\theta_i$ and $r_k$ using \eqref{cosine law} as follows:
\begin{equation}
\begin{split}
\tilde{r}_{i,k}^2=r_{k}^2+D^2-2{r_k}D \;\mathrm{cos}{\theta_i} \:\:\: \forall r_k , \forall \theta_i,
\end{split}
\end{equation}
where, $\tilde{r}_{i,k}$ denotes the interfering distance from a specified polar coordinate $( r_k,\theta_i)$ in the interfering cell to the BS of interest located at a distance $D$ (see \figref{Cap0}).
In addition, it is worth to mention that $\tilde{r}_{i,k}$ are the points at which $P(\tilde{r}_{\mathrm{sel}}=\tilde{r}_{i,k})$ can be defined using \eqref{IntPMF} as follows:
\begin{equation}
P(\tilde{r}_{\mathrm{sel}}=\tilde{r}_{i,k})\\=\frac{P(r_{\mathrm{sel}}=r_k)}{\mathcal{I}}.
\end{equation}
The two dimensional data set of $\tilde{r}_{\mathrm{sel}}$, at which $P(\tilde{r}_{\mathrm{sel}}=\tilde{r}_{i,k})$ is defined, can then be grouped into $M$ segments of any arbitrary width $\Delta$. This can be done by dividing the distance between $D-R$ and $D+R$ into $M$ equal segments of width $\Delta$\textcolor{black}{\footnote{\textcolor{black}{Note that $\Delta$ represents the uniform segments of interfering cell whereas $\Delta_k$ represents the non-uniform circular regions within a given cell.}}} and mapping $\tilde{r}_{i,k}$ accordingly. Clearly, by adding all the probabilities for which $\tilde{r}_{\mathrm{sel}}$ lies in the $m^{\mathrm{th}}$ segment we get the probability of $\tilde{r}_{\mathrm{sel}}=\tilde{r}_m$:
\begin{equation}
\label{step}
P(\tilde{r}_{\mathrm{sel}}=\tilde{r}_m)= \sum_{\tilde{r}_{i,k} \in[\tilde{r}_m - \frac{\Delta}{2}, {\tilde{r}_m+\frac{\Delta}{2}}]}
P(\tilde{r}_{\mathrm{sel}}=\tilde{r}_{i,k}),
\end{equation}
where $\tilde{r}_m$ denotes any discrete value that the RV $\tilde{r}_{\mathrm{sel}}$ can take. Recall $X=\bar{K}\tilde{r}^{-\beta} \chi$, therefore the PDF of $X$ conditioned on $\tilde{r}_{\mathrm{sel}}$ can be determined by RV transformation as:
\begin{equation}
f_{X|\tilde{r}_{\mathrm{sel}}}=\frac{f_\chi(x\tilde{r}_{\mathrm{sel}}^{\beta}\bar{K}^{-1})}{\bar{K}\tilde{r}_{\mathrm{sel}}^{-\beta}}.
\end{equation}
Averaging over the PMF of $\tilde{r}_{\mathrm{sel}}$, the distribution of the ICI, $f_X(x)$, from any cell $l$ can be given as:
\begin{equation}
\label{ICI}
f_{X}(x)=\sum_{\tilde{r}_m=\tilde{r}_1}^{\tilde{r}_M} \frac{f_\chi(x\tilde{r}_m^{\beta}\bar{K}^{-1})}{\bar{K}\tilde{r}_m^{-\beta}}P(\tilde{r}_{\mathrm{sel}}=\tilde{r}_m).
\end{equation}
It is important to emphasize that the derivation of the distribution of ICI is based on the scheduling decisions of interfering cells at a given time slot. Therefore, the parameter $\tilde{r}_m$ of the ICI distribution varies from one time slot to the other for the location based round robin and greedy round robin schemes.

\section{MGF of the Cumulative Intercell Interference}
Computing the distribution of the cumulative ICI $Y$  requires the convolution of the PDF of $L$ RVs $X_l$, $\forall l=1,2,\cdots L$, which is a tedious task in practice. To avoid convolutions, we utilize an MGF based approach to derive the expression for the MGF of the cumulative ICI $Y$. Considering same scheduling scheme deployed in each cell, the interference experienced from each cell is i.i.d and therefore the MGF of the cumulative interference  can be calculated as follows:
\begin{equation}
\mathcal{M}_{Y}(t)=\prod_{l=1}^{L} \mathcal{M}_{X_l}(t)=\left(\mathcal{M}_{X}(t)\right)^L=\left(\mathbb{E}[e^{tx}]\right)^L.
\end{equation}
The derivation can be extended in a straightforward manner to the case where different cells deploy different scheduling schemes\footnote{The only required change is that in this case, the  MGF of  ICI received from each cell will not be i.i.d, i.e.,  $\mathcal{M}_Y(t)=\prod_{l=1}^L \mathcal{M}_{X_l}(t)$ and the MGF of each cell $l$ can be characterized with its corresponding scheduling scheme which depends on $P(r_{\mathrm{sel}}=r_k)$.}.
Looking at the structure of \eqref{ICI}, we can write:
\begin{equation}\label{cumMGF}
\begin{split}
\mathcal{M}_{X}(t)=\int_0^\infty e^{tx}f_{X}(x) dx=\sum_{\tilde{r}_m=\tilde{r}_1}^{\tilde{r}_M}
\frac{P(\tilde{r}_{\mathrm{sel}}=\tilde{r}_m)}{\bar{K}\tilde{r}_m^{-\beta}} \times\\ \int_0^\infty  e^{tx}{f_\chi(x\tilde{r}_m^{\beta}\bar{K}^{-1})}dx.
\end{split}
\end{equation}
The derived expression is generic and applies to any composite fading distribution. Next, we will present explicit MGF expressions for three typically used practical fading models.\\
{\it \bf Special Case 1: Rayleigh Fading -$\zeta,  \chi \sim \mathrm{Exp}(\lambda)$:}
In this case, the small scale fading coefficient on a given subcarrier is considered to be Rayleigh distributed whereas the effect of shadowing is not considered. The distribution of interference considering a single interfering cell can then be derived as:
\begin{equation}
\label{HEdist}
f_{X}(x) = \sum_{\tilde{r}_m=\tilde{r}_1}^{\tilde{r}_M}\bar{K}^{-1}\lambda\tilde{r}_m^{\beta} e^{-\lambda \tilde{r}_m^\beta \bar{K}^{-1} x} P(\tilde{r}_{\mathrm{sel}}=\tilde{r}_m).
\end{equation}
Note that \eqref{HEdist} is a Hyper-Exponential distribution with parameter $\bar{K}^{-1} \lambda \tilde{r}_m^\beta$. Thus, using the MGF of the Hyper-Exponential distribution, $\mathcal{M}_{Y}(t)$ can be derived as follows:
\begin{equation}
\label{finalmgfHE}
\mathcal{M}_{X}(t)=\left(\sum_{\tilde{r}_m=\tilde{r}_1}^{\tilde{r}_M}\frac{\bar{K}^{-1}\lambda \tilde{r}_m^\beta}{\bar{K}^{-1}\lambda \tilde{r}_m^\beta-t}P(\tilde{r}_{\mathrm{sel}}=\tilde{r}_m)\right)^L.
\end{equation}
{\it  \bf Special Case 2: Generalized-$\mathcal{K}$ Composite Fading -{$\zeta,  \chi \sim \mathcal{K}_G(m_s,m_c, \Omega)$:}}
In wireless channels, shadowing and fading across the channel between a user and BS can be jointly modeled by a composite fading distribution.
A closed form composite fading model, namely Generalized-$\mathcal{K}$ also referred to as Gamma-Gamma distribution, has been recently introduced in \cite{KG} which is general enough to model well-known shadowing and fading distributions such as log-normal, Nakagami-\emph{m}, Rayleigh etc.  Using \eqref{ICI}, $f_{X}(x)$ in this case can be derived as follows:
\begin{equation}
\begin{split}
f_{X}(x) = \sum_{\tilde{r}_m=\tilde{r}_1}^{\tilde{r}_M}  \frac{2(x{\tilde{r}_m}^{\beta}\bar{K}^{-1})^{\frac{m_c+m_s-2}{2}}}{\bar{K}\tilde{r}_m^{-\beta}\Gamma(m_c)\Gamma(m_s)} \left(\frac{b}{2}\right)^{m_c+m_s} \times \\ \mathbb{K}_{m_s-m_c}\left(b\sqrt{x\tilde{r}_m^\beta \bar{K}^{-1}}\right)P(\tilde{r}_{sel}=\tilde{r}_m),
\end{split}
\end{equation}
where, $\mathbb{K}_v(.)$ denotes the modified Bessel function of second kind with order $v$, $b=2\sqrt{\frac{m_c m_s}{\Omega}}$. Performing some algebraic manipulations and using  \cite[Eq. 6.643/3]{book}, the expression for $\mathcal{M}_{X}(t)$ can be derived as:
\begin{equation}
\begin{split}
\mathcal{M}_{X}(t)=\sum_{\tilde{r}_m=\tilde{r}_1}^{\tilde{r}_M}
P(\tilde{r}_{sel}={\tilde{r}_m}^\beta) e^{\frac{b^2\tilde{r}_m^{\beta}}{8\bar{K}t}} \left(\frac{-b^2\tilde{r}_m^{\beta}}{4\bar{K}t}\right)^{(\frac{m_s+m_c-1}{2})} \times \\
\mathbb{W}_{\frac{1-m_c-m_s}{2},\frac{m_c-m_s}{2}}\left(-\frac{\tilde{r}_m^{\beta}b^2}{4\bar{K}t}\right),
\end{split}
\end{equation}
where, $\mathbb{W}$ denotes the Whittaker function.
Finally, $\mathcal{M}_{Y}(t)$ can be written as follows:
\begin{equation}
\label{finalmgfHG1}
\begin{split}
\mathcal{M}_{Y}(t)=\left(\sum_{\tilde{r}_m=\tilde{r}_1}^{\tilde{r}_M}
P(\tilde{r}_{sel}={\tilde{r}_m}^\beta)
\left(-\frac{\tilde{r}_m^{\beta}b^2}{4\bar{K}t}\right)
e^{\frac{b^2\tilde{r}_m^{\beta}}{8\bar{K}t}} \times \right.
\\\left.
\mathbb{W}_{\frac{1-m_c-m_s}{2},\frac{m_c-m_s}{2}}
\left(\frac{-b^2\tilde{r}_m^{\beta}}{4\bar{K}t}\right)^{(\frac{m_s+m_c-1}{2})}\right)^L.
\end{split}
\end{equation}
Note that integrating the CDF and MGF of Generalized-$\mathcal{K}$ RV which involves Meijer-G and Whittaker functions, respectively, in \texttt{MATHEMATICA} and \texttt{MAPLE} can be a bit more time consuming. In this context,  recently an  approximation of the Generalized-$\mathcal{K}$ RV has been proposed in \cite{KG} which is discussed below.\\
{\it \bf Special Case 3: Gamma Composite Fading -{$\zeta, \chi \sim \mathrm{Gamma}(m_s,m_c)$}}:
In \cite{KG}, the authors proposed an accurate approximation of the Generalized-$\mathcal{K}$ RV by the more tractable Gamma distribution  using  moment matching method.
The approximation  provides a simplifying model for the composite fading in wireless communication systems.
Using \eqref{ICI}, $f_{X}(x)$ can be written in this case as:
\begin{equation}
f_{X}(x) = \sum_{\tilde{r}_m=\tilde{r}_1}^{\tilde{r}_M} \frac{e^{-\frac{x \tilde{r}_m^\beta \bar{K}^{-1}}{m_c}}(x \tilde{r}_m^\beta \bar{K}^{-1})^{m_s-1}}{\bar{K} \tilde{r}_m^{-\beta} \Gamma(m_s) m_c^{m_s}} P(\tilde{r}_{sel}=\tilde{r}_m).
\end{equation}
Performing some algebraic manipulations and letting $y=x(\frac{\tilde{r}_m^{\beta}}{m_c}-t)$, $\mathcal{M}_{X}(t)$ can be derived as follows:
\begin{equation}
\begin{split}
\mathcal{M}_{X}(t)=\sum_{\tilde{r}_m=\tilde{r}_1}^{\tilde{r}_M}  \frac{P(\tilde{r}_{sel}=\tilde{r}_m)(\bar{K}^{-1}\tilde{r}_m^{\beta})^{m_s}}{\Gamma(m_s) 
\left( {\tilde{r}_m^{\beta}\bar{K}^{-1}}-{m_c}t\right)^{m_s}}
\int_0^\infty  e^{-y}y^{m_s-1}dy\\=\sum_{\tilde{r}_m=
\tilde{r}_1}^{\tilde{r}_M}  P(\tilde{r}_{sel}=\tilde{r}_m)\left(\frac{\bar{K}^{-1}\tilde{r}_m^{\beta}}{\bar{K}^{-1}\tilde{r}_m^{\beta} -m_ct}\right)^{m_s}.
\end{split}
\end{equation}
Finally, $\mathcal{M}_{Y}(t)$ can be written as follows:
\begin{equation}
\label{finalmgfHG2}
\mathcal{M}_{Y}(t)=\left(\sum_{\tilde{r}_m=\tilde{r}_1}^{\tilde{r}_M} P(\tilde{r}_{\mathrm{sel}}=\tilde{r}_m)\left(\frac{\bar{K}^{-1}\tilde{r}_m^{\beta}}{\bar{K}^{-1}\tilde{r}_m^{\beta} -m_ct}\right)^{m_s}\right)^L.
\end{equation}

\section{Evaluation of Important Network Performance Metrics}
In this section, we demonstrate the significance of the derived  MGF expressions of the cumulative ICI  in quantifying important network performance metrics such as the outage probability $P_{\mathrm{out}}$, ergodic capacity $\mathcal{C}$ and average fairness $\mathcal{F}$ among users numerically.\\
\\
\noindent
{\it \bf Evaluation of Outage Probability:}
The outage probability is typically defined as the probability of the instantaneous interference-to-signal-ratio to exceed a certain threshold. In order to evaluate $P_{\mathrm{out}}$ numerically, we use the MGF based technique introduced in  \cite{outage} for interference limited systems.
Firstly, we define a new RV, ${Z}=q\sum_{l=1}^{L} X_l-X_0=qY-X_0$, where $q$ is the outage threshold and $X_0$ is the corresponding signal power  of the scheduled user  in the central cell. An outage event occurs when $p({Z} \geq 0)$, i.e., when the interference exceeds the corresponding signal power. This decision problem is solved in \cite{outage} by combining the characteristic function of $Z$ and residue theorem. The characteristic function of ${Z}$ is defined as $\phi_{Z}(j\omega)=\mathbb{E}[e^{j{Z}\omega}]$. Considering interference $Y$ and signal power $X_0$ to be independent, the expression for $\phi_{Z}(j\omega)$ can be given as,
$
\phi_{Z}(j\omega)=\phi_{Y}(jq\omega)\phi_{X_0}(-j\omega)
$; where $\phi_Y(q j\omega)$ can be given by \eqref{finalmgfHE}, \eqref{finalmgfHG1}, and \eqref{finalmgfHG2} for different fading models. In general, the characteristic function of $X_0$ can be calculated as \cite{yilmaz}:
\begin{equation}\label{xo}
\begin{split}
\phi_{X_0}(\omega)=\mathbb{E}(e^{j\omega x_0})=\int_0^\infty e^{j \omega x_0} f_{X_0}(x_0) d x_0
\\=j \omega\int_0^\infty e^{j \omega x_0} F_{X_0}(x_0) d x_0,
\end{split}
\end{equation}
where the detailed expressions of $F_{X_0}(x_0)$ for all opportunistic scheduling schemes are provided in \cite{hina}.
Compact closed form expressions  of $\phi_{X_0}(j \omega)$ are available in the literature for e.g. \cite[Eq. 19]{romero}. 
For non-opportunistic scheduling schemes
$
\phi_{X_0}(j\omega)= \sum_{{r}_k={r}_1}^{{r}_K}  \phi_{\zeta|r_k}(j \omega) P(r_{\mathrm{sel}}={r}_k)
$,
where $\phi_{\zeta|r_k}(j \omega)$ is the characteristic function of $\zeta$ in ring $k$.  The outage probability can then be computed by using the  classical lemma introduced in \cite{outage} as follows:
\begin{equation}
\label{QTZhang}
P_{\mathrm{out}}=\frac{1}{2} +\frac{1}{\pi }\int_{0}^\infty \mathrm{Im}\left(\frac{\phi_Z(\omega)}{\omega}\right) d\omega,
\end{equation}
where $\mathrm{Im}(\phi_Z(\omega))$ denotes the imaginary part of $\phi_Z(\omega)$. 
Using \eqref{QTZhang}, the outage probability can be evaluated using any standard mathematical software packages such as \texttt{MATHEMATICA} and \texttt{MAPLE}. 
\\
\\
\noindent
{\it \bf Evaluation of Ergodic Network Capacity:}
Another important performance evaluation parameter is the network ergodic capacity $\mathcal{C}$, i.e.,
$
\mathcal{C}=\mathbb{E}\left[\mathrm{log_2}\left(1+\frac{X_0}{\sum_{l=1}^L{X}_l+\sigma^2}\right)\right]
$
Usually, the computation of $\mathcal{C}$ requires $(L+1)$-fold numerical integrations. To avoid this, we utilize the efficient lemma derived in \cite{hamdi} with a slight modification to take thermal noise into account as follows:
\begin{equation}
\begin{split}\label{lemma}
&
\mathbb{E}\left[\mathrm{ln}\left(1+\frac{X_0}{\sum_{l=1}^L{X}_l+\sigma^2}\right)\right]\\&=\int_0^\infty \frac{\mathcal{M}_Y(t)-\mathcal{M}_{X_0,Y}(t)}{t} e^{-\sigma^2 t} dt,
\end{split}
\end{equation}
where, $\mathcal{M}_Y(t)=\mathbb{E}[e^{-t\sum_{l=1}^{L} X_l}]$ and $\mathcal{M}_{X_0,Y}(t)=\mathbb{E}[e^{-t(X_0+\sum_{l=1}^{L} X_l)}]=\mathbb{E}[e^{-t(X_0+Y)}]$. Note that this is the definition of MGF as defined in \cite{hamdi} which is not the same as our definition. Thus,
we can use $\mathcal{M}_Y(t)$ from \eqref{finalmgfHE}, \eqref{finalmgfHG1} and \eqref{finalmgfHG2} directly with a sign change of $jw$. 
In addition, \eqref{lemma} can be solved efficiently by expressing it in terms of the weights and abscissas of a Laguerre orthogonal polynomial \cite{hamdi} as
$
\mathbb{E}\left[\mathrm{ln}\left(1+\frac{X_0}{\sum_{l=1}^L{X}_l+1}\right)\right]=
\sum_{\epsilon =1}^E \alpha_\epsilon \frac{\mathcal{M}_Y(\xi_\epsilon)-\mathcal{M}_{X_0,Y}(\xi_\epsilon)}{\xi_\epsilon} +R_E
$,
where $\xi_\epsilon$ and $\alpha_\epsilon$ are the sample points and the weight factors of the Laguerre polynomial, tabulated in \cite{book}, and $R_E$ is the remainder. The MGF of $X_0$ can be calculated as explained in \eqref{xo}.
\\
{\it \bf Evaluation of Average Fairness:}
In order to quantify the degree of fairness among different scheduling schemes, we use the notion developed in \cite{elliott}. The average fairness of a scheduling scheme with $U$ users can be given as,
$
\mathcal{F}=-\sum_{i=1}^{U} p_i\frac{ \mathrm{log}_{10} p_i}{\mathrm{log}_{10}U}
$,
where $p_i$ is the proportion of resources allocated to a user $i$ or the access probability of user $i$. A system is strictly fair if each user has equal probability to access the channel and in such case the average fairness becomes one. The other extreme occurs when the channel access is dominated by a single user; in such case, the average fairness reduces to zero. The average fairness can be easily computed using our derived results as 
\begin{equation}
\mathcal{F}=-\sum_{k=1}^{K} P(r_{\mathrm{sel}}=r_k) \frac{ \mathrm{log}_{10} P(r_{\mathrm{sel}}=r_k)- \mathrm{log}_{10} {u_k}}{\mathrm{log}_{10}U},
\end{equation}
where $u_k$ denotes the number of users in a ring $k$.
\section{Numerical and Simulation Results}
In this section, we first define the system parameters and discuss the Monte-Carlo simulation setup which is required to demonstrate the accuracy of the derived expressions.  
We then address some important insights and study the performance trends of the  different scheduling schemes.
\subsection{ Parameter Settings and Simulation Setup:}
The radius $R$ of the cell is set to 500m and the cell is decomposed into non-uniform circular regions of width $\Delta_k$. The path loss variation within each circular region is set to $\kappa=2$dB.

For each Monte-Carlo trial, we generate $U$ uniformly distributed users in a circular cell of radius $R$. Each user has instantaneous SNR given by \eqref{1} and \textcolor{black}{short term average selected SNR based on its ring location $(\bar{\gamma}_k)$}.   We allocate a user with maximum instantaneous SNR in the greedy scheme, a user with maximum normalized SNR in the proportional fair scheme and a user arbitrarily for the round robin scheme.  For location based round robin  we select a user randomly from the $w^{th}$ ring in a time slot $w$ whereas we select a user with maximum SNR  without considering  the users of the previously allocated $w-1$ rings for the greedy round robin scheme. 
Next, we calculate the distances of the selected users, i.e., $r_{\mathrm{sel}}$, from their serving BS and compute the distances to the BS of interest, i.e., $\tilde{r}_{\mathrm{sel}}$ for all scheduling schemes. 
The process repeats for large number of Monte-Carlo  simulations. The  distance data is then analyzed by creating a histogram of non-uniform and uniform bin width \textcolor{black}{in \figref{Cap1} and  \figref{Cap2}, respectively}. 

\subsection{ Results and Discussions:}
\begin{figure}[t]
  \centering
  \includegraphics[totalheight=3in,width=3.65in]{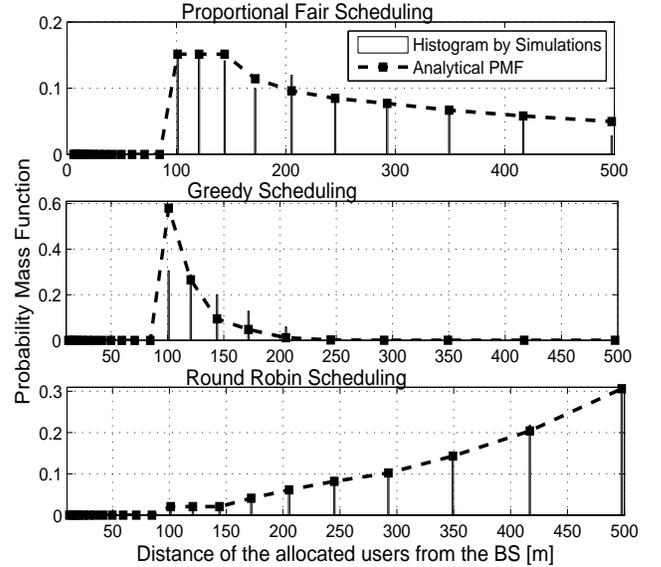}
\caption{PMF of the distance of the allocated users in a given cell (i.e., PMF of $r_{\mathrm{sel}}$) for proportional fair, greedy, and round robin scheduling schemes with path loss exponent $\beta=2.6$, $U=50$, $C$=60 dB, $P_{\mathrm{max}}$=1W, $\sigma^2$=-174 dBm/Hz, and Number of Monte-Carlo simulations =100,000.}
\vspace{-3mm}
\label{Cap1}
\end{figure}

\figref{Cap1} depicts the PMF of the location of the scheduled user in a given cell based on the proportional fair, greedy and round robin scheduling schemes. Since the proportional fair scheme exhibits some fairness among users in a cell, the PMF of the allocated user locations is expected to be more flat compared to the greedy scheme. Since the cell edge has more users due to the large area and each user has equal probability to be allocated on a given subcarrier, therefore the round robin scheme exhibits high probability at the cell-edge. In order to get an integer number of users within a ring, we perform rounding in the analysis, i.e., we consider zero active users in the rings where $u_k \leq 0.5$. In Monte-Carlo simulations, we consider the probability of allocating a user in these rings to be zero which can also be verified from \figref{Cap1}.  The impact of rounding  can be observed more in the  cell center  than at the cell-edge, due to narrow circular regions  and  in turn low number of users in each ring in the cell center. This rounding effect is therefore more visible in the PMF of the greedy scheme which exist mainly in the cell center  as observed in Fig.~2(b).  However, as this mismatch lies near the cell center,  it does not have much impact on the  derived performance results  which include distribution of the ICI, ergodic capacity and outage probability.

It is important to note that the numerical results for the derived PMF in \figref{Cap1} nearly coincide with the exhaustive Monte-Carlo simulation results with a small number of rings $K=10$ and  $\kappa=2$dB. 
Moreover, it can also be noticed that the width of the circular regions tend to increase from cell center to cell edge which is due to the exponentially decaying path loss as mentioned in Section II. The number of required rings is expected to decrease by reducing $\beta$ and  increasing the amount of power decay within each circular region and vice versa.
Another important point to explain with reference to \figref{Cap1} is that with the increase in the number of competing users on a given subcarrier, the PMF of opportunistic scheduling schemes will  get more skewed toward the cell center which is due to the fact that the higher the  number of users in the cell center, the higher is the probability of allocating a subcarrier in the cell center.

\begin{figure}[t]
 \centering
  \includegraphics[totalheight=3in,width=3.65in]{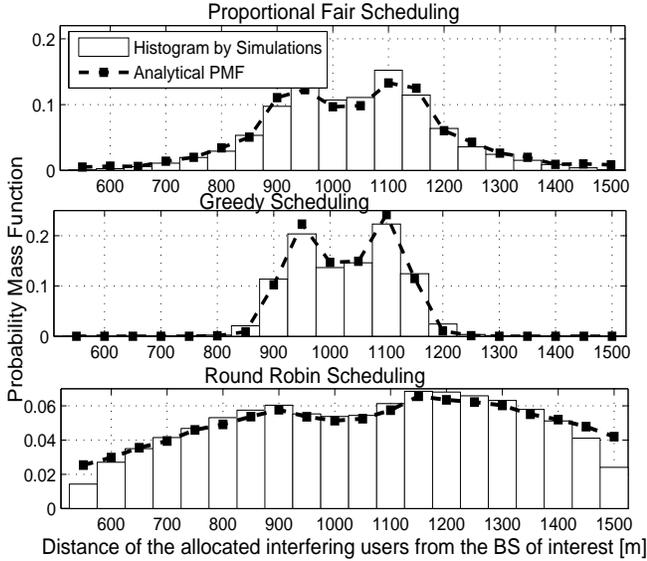}
   \caption{PMF of the distance at which the users in the interfering cells are allocated (i.e., PMF of $\tilde{r}_{\mathrm{sel}}$) for proportional fair, greedy and round robin scheduling schemes with path loss exponent $\beta=2.6$, $U=50$, $\mathcal{I}=180$, $\chi \sim \mathrm{Gamma}(3/2,2/3)$, $C$=60 dB, $P_{\mathrm{max}}$=1W, $\sigma^2$=-174 dBm/Hz, $\Delta$=50 m, and Number of Monte-Carlo simulations =100,000.}
\vspace{-3mm}
\label{Cap2}
\end{figure}

In \figref{Cap2}, the PMF of the distance between the allocated user in interfering cell $l$ and the BS of interest, i.e., $P(\tilde{r}_{\mathrm{sel}}=\tilde{r}_m)$, is presented. Numerical results are found to be in close agreement with the Monte-Carlo simulation results \textcolor{black}{and $M=20$, i.e., $\Delta$=50 m}. 
For the opportunistic scheduling schemes, it is likely that a user close to its serving BS can get a subcarrier, thus, the PMF of the distance of allocated interfering users is expected to have high density in the middle. However, the  slight descend in the central region in \figref{Cap2} is due to ignoring users that lie within the rings where the average number of users is less than half.
Moreover, we can observe that the round robin scheduler is highly vulnerable to interference compared to the other schemes as high interference is expected to come from the cell edge users in the interfering cells. On the other hand, the greedy scheduler is expected to have allocations near the cell center and, thus, leads to less interference from neighboring cells. The proportional fair scheme lies in between the two extremes. \textcolor{black}{Note that as $M$ increases the analytical  accuracy is expected to increase even further with a trade-off  in terms of computational complexity.  An acceptable range of $\Delta$ with negligible complexity varies from 1 to 100 m}.

\begin{figure}
  \centering
 \includegraphics[totalheight=3in,width=3.65in]{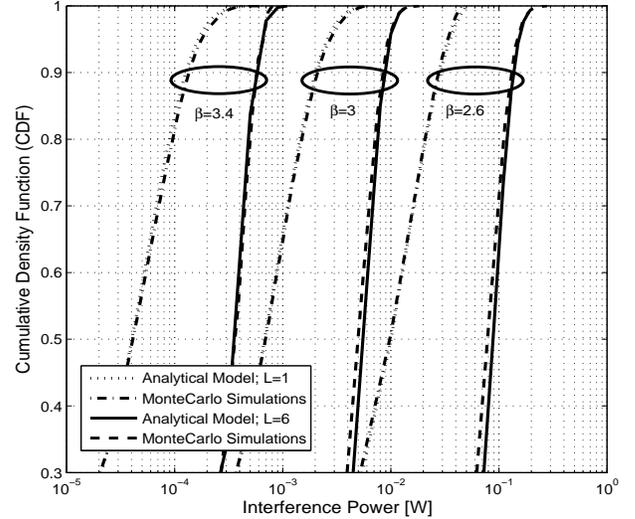}
   \caption{Impact of different number of interferers $L$ and various path loss exponents $(\beta)$ on the CDF of cumulative ICI  considering greedy scheduling scheme, $U=50$, $\mathcal{I}=180$, $\chi \sim \mathrm{Gamma}(3/2,2/3)$, $C$=60 dB, $P_{\mathrm{max}}$=1W, $\sigma^2$=-174 dBm/Hz,  and Number of Monte-Carlo simulations =100,000.}
\vspace{-3mm}
\label{Cap3}
\end{figure}

\figref{Cap3} illustrates the CDF of the ICI considering different number of interfering cells and path loss exponents $\beta$ for the greedy scheduling scheme. With the increase in the number of interferers, the interference level increases.  Moreover, as $\beta$ increases, the signal degrades rapidly and thus interference level is reduced considerably.
At this point, it is important to mention  that in this paper we derive and utilize the  MGF of the cumulative ICI rather than the  CDF of the cumulative ICI in order to evaluate important network performance metrics. Therefore, the analytical part of the provided figure of the CDF of the cumulative ICI is plotted   using a  technique mentioned in \cite{sir} to convert MGF into CDF numerically.

\begin{figure}[t]
  \centering
  \includegraphics[totalheight=3in,width=3.65in]{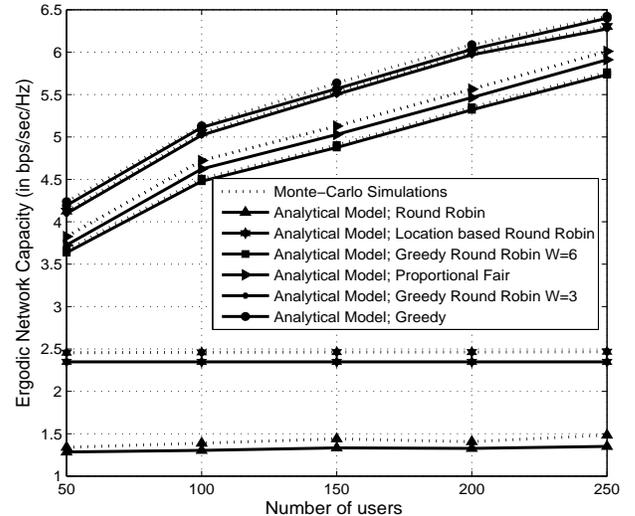}
   \caption{Network ergodic capacity  for different number of users considering different scheduling schemes with $\beta=2.6$, $\chi \sim \mathrm{Gamma}(3/2,2/3)$, $C$=60 dB, $P_{\mathrm{max}}$=1W, and $\sigma^2$=-174 dBm/Hz.}
\vspace{-3mm}
\label{Cap4}
\end{figure}

\figref{Cap4} investigates the effect of increasing the number of competing users on a given sub-carrier considering all scheduling schemes. It can be observed that the increase in the number of users enhances the performance of the opportunistic scheduling schemes due to additional multiuser diversity gains.  
The greedy scheme achieves the best performance whereas the round robin scheme achieves the worst performance. As expected, the proportional fair scheme lies in between the two extremes. The average capacity of location based round robin over $W=K$ time slots has been shown  to be better than the conventional round robin scheme. The average capacity results of the greedy round robin scheme is presented for $W=3$ and $W=6$. Clearly, for $W=1$, the scheme is equivalent to the greedy scheme; however, with the increase of time slots, performance degradation takes place due to the reduction of  multiuser diversity caused by ignoring the users from  previously allocated rings. 

\begin{figure}[t]
  \centering
  \includegraphics[totalheight=3in,width=3.65in]{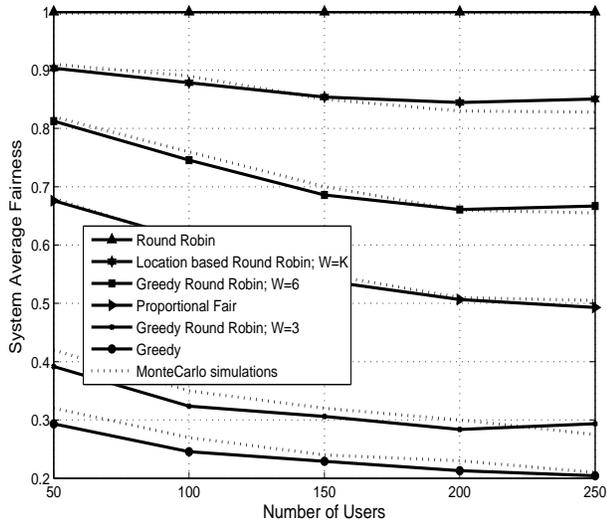}
   \caption{Average system fairness for different number of users considering different scheduling schemes with $\beta=2.6$, $C$=60 dB, $P_{\mathrm{max}}$=1~W, and $\sigma^2$=-174 dBm/Hz.}
\vspace{-3mm}
\label{Cap5}
\end{figure}

\figref{Cap5} quantifies the average resource fairness of all  presented scheduling schemes. As expected, round robin is a strictly fair scheme. The proportional fair scheme possesses the ability to enhance the network throughput compared to round robin scheduling while providing a high degree of fairness. The greedy scheme is observed to be the most unfair scheme. Considering $K$ time slots, the average fairness of the location based round robin scheme is investigated and found to be very close to the round robin scheme, however, with degradation in performance as can be observed in \figref{Cap4}. For the greedy round robin scheme, we plotted the fairness metric considering $W=3$ and $W=6$; it is shown that as the number of time slots increases, the fairness improves with a trade-off price in terms of ergodic capacity.

\begin{figure}[t]
  \centering
 \includegraphics[totalheight=3in,width=3.65in]{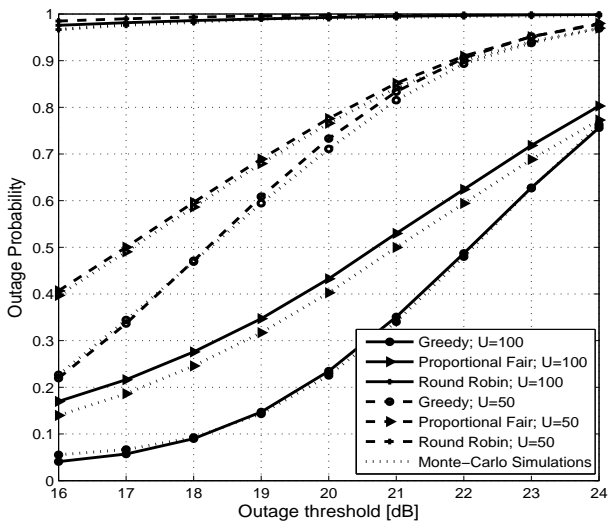}
  \caption{Impact of different scheduling schemes on the network outage probability for different number of users and various  outage thresholds, $\beta$=2.6, $\chi \sim \mathrm{Gamma}(3/2,2/3)$, $C$=60 dB, $P_{\mathrm{max}}$=1W, and $\sigma^2$=-174 dBm/Hz.}
\vspace{-3mm}
\label{Cap6}
\end{figure}

\begin{figure}[t]
  \centering
  \includegraphics[totalheight=3in,width=3.65in]{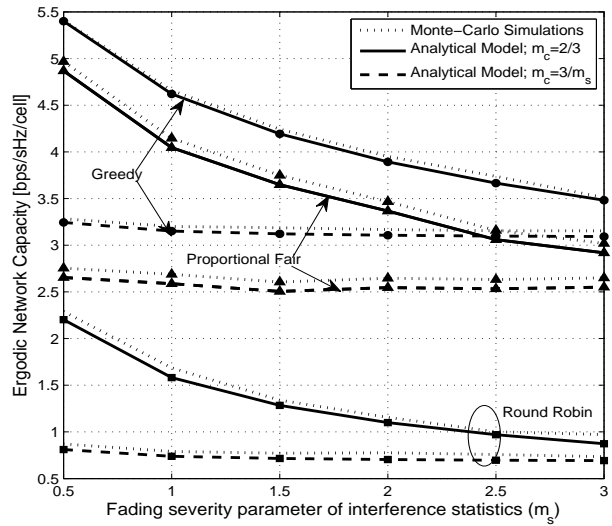}
   \caption{Ergodic network capacity  of different scheduling schemes as a function of  different  parameters of the interference statistics $\chi \sim \mathrm{Gamma} (m_s,m_c)$  with $U =50$ users, $\beta$=2.6, and $\sigma^2$=-174 dBm/Hz.}
\vspace{-3mm}
\label{Cap7}
\end{figure}
In \figref{Cap6}, we evaluate the network outage probability as a function of the outage threshold; $q=(Z+X_0)/Y$ for (i) $U=50$ users; (ii) $U=100$  users. The numerical and simulation results are nearly identical for most cases.
 The higher the outage threshold for a given  signal and interference power, the greater outage  is expected. Moreover, for larger number of users the outage probability is observed to reduce for all opportunistic scheduling schemes except the round robin scheme. Since increasing the number of users on a given subcarrier in non-opportunistic schemes does not directly affect the access probability of a ring $k$, therefore its impact on the ICI is almost negligible. This fact can also be verified from \figref{Cap4}.
Finally, in \figref{Cap7}, we evaluate the network ergodic capacity as a function of the fading severity parameter and average power of Gamma fading interference channels for different scheduling schemes. Firstly, it can be observed that  increasing the average power of the interference channel which is given by $\Omega=m_c m_s$ for a given fading severity parameter $m_s$, the capacity degrades significantly for all schemes. Moreover, it is also shown that increasing the fading severity $m_s$ while keeping the average power $\Omega=3$ fixed has minimal impact on the system capacity. Therefore, the lower average power of interference channel $\Omega$, the better is the overall system performance.

The small gap between the analytical and simulation results is mainly due to assuming in the analytical derivations that users located within a ring are at the boundary of the ring. \textcolor{black}{ This gap can be further reduced in \figref{Cap1} by  increasing the accuracy of the approximation, i.e., by reducing the path loss decay $\kappa$ within each circular region which in turn increases the number of rings. An acceptable range of $\kappa$ with reasonable complexity varies from 0.5 to 3 dB such that $r \leq R$.}

\section{Conclusion}
We proposed a novel approach to model the uplink ICI considering various scheduling schemes and composite fading channel models. The proposed approach is not  dependent on a particular shadowing and fading statistics, hence,  extensions to different models is possible.  The provided \textcolor{black} {numerical results} help in gaining insights into the behavior of ICI considering different scheduling schemes and composite fading models. Moreover, they provide quantitative assessment of the relative performance of various scheduling schemes which is important for  network design and assessment. The proposed approach can be extended to typical power control schemes and downlink scenarios as discussed in \cite{hina}.

\section{Acknowledgment}
This work  was made possible by NPRP grant 4-353-2-130 from the Qatar National Research Fund (a member of The Qatar Foundation). The statements made herein are solely the responsibility of the authors. 
\bibliography{IEEEfull,references}
\bibliographystyle{IEEEtran}

\begin{IEEEbiography}[{\includegraphics[width=1in,height=1.25in,clip,keepaspectratio]{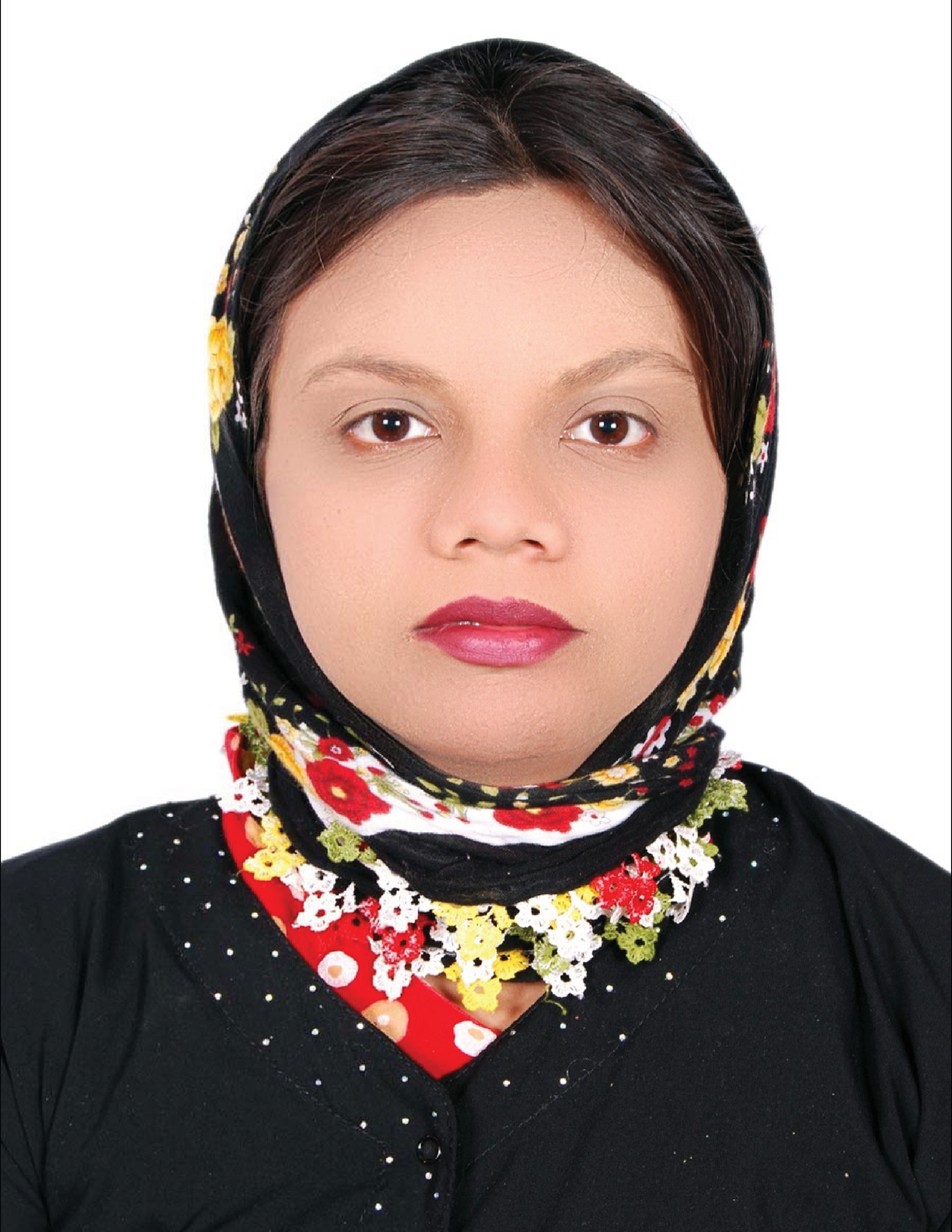}}]
{Hina Tabassum} 
received her Bachelors degree in Electronics from the N.E.D University of Engineering and Technology (NEDUET), Karachi, Pakistan, in 2004. She received during her undergraduate studies the Gold medal from NEDUET and from SIEMENS for securing the first position among all engineering universities of Karachi. She then worked as lecturer in NEDUET for two years. In September 2005, she joined the Pakistan Space and Upper Atmosphere Research Commission (SUPARCO), Karachi, Pakistan and received there the best performance award in 2009. She also completed her Masters in Communications Engineering from NEDUET in 2009. In January 2010, she joined the Computer, Electrical, and Mathematical Sciences \& Engineering Division at King Abdullah University of Science and Technology (KAUST), Thuwal, Makkah Province, Saudi Arabia, where she is currently a Ph.D. candidate. Her research interests include wireless communications with focus on  interference modeling, radio resource allocation, and optimization in heterogeneous networks.
\end{IEEEbiography}

\begin{IEEEbiography}[{\includegraphics[width=1in,height=1.25in,clip,keepaspectratio=false]{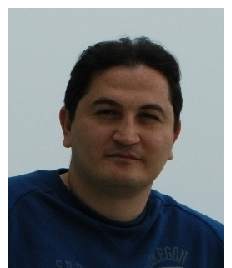}}]{Ferkan Yilmaz} (M'10) was born in Malatya,
Turkey. He received the B.Sc. degree in electronics and
communications (with the first and the second highest honors in the
electronics \& communications engineering department and the
university, respectively) from Yildiz Technical University (YTU),
Istanbul, Turkey, in 1997, and the M.Sc. degree in electronics and
communications from Istanbul Technical University, Istanbul, Turkey,
in 2002. He received his PhD degree from the Telecommunications
Branch at Gebze Institute of Technology (GYTE), Turkey, in January
2009, where received the award for the best PhD thesis in October
2009. From 1998 to 2003, he worked for the National Research
Institute of Electronics and Cryptology, Tubitak, Turkey. From 2004
to 2008, he worked for Vodafone Technology in Turkey as a senior
telecommunications researcher. He was a research associate for the
Texas A\&M University at Qatar in 2008, and is currently a post
doctoral fellow at King Abdullah University of Science and
Technology (KAUST) where he has worked since August 2009. His
research interests include digital signal processing in
communications, signal propagation aspects and diversity reception
techniques in wireless mobile radio systems, spread-spectrum
techniques, cooperative/collaborative communications, and multihop
communications. Moreover, he is interested in wireless nano-circuit
clouds/fractals and wireless fractional communications techniques
with emphasis on combinatorics, special functions, and linear
transformations.
\end{IEEEbiography}

\begin{IEEEbiography}[{\includegraphics[width=1in,height=1.25in,clip,keepaspectratio]{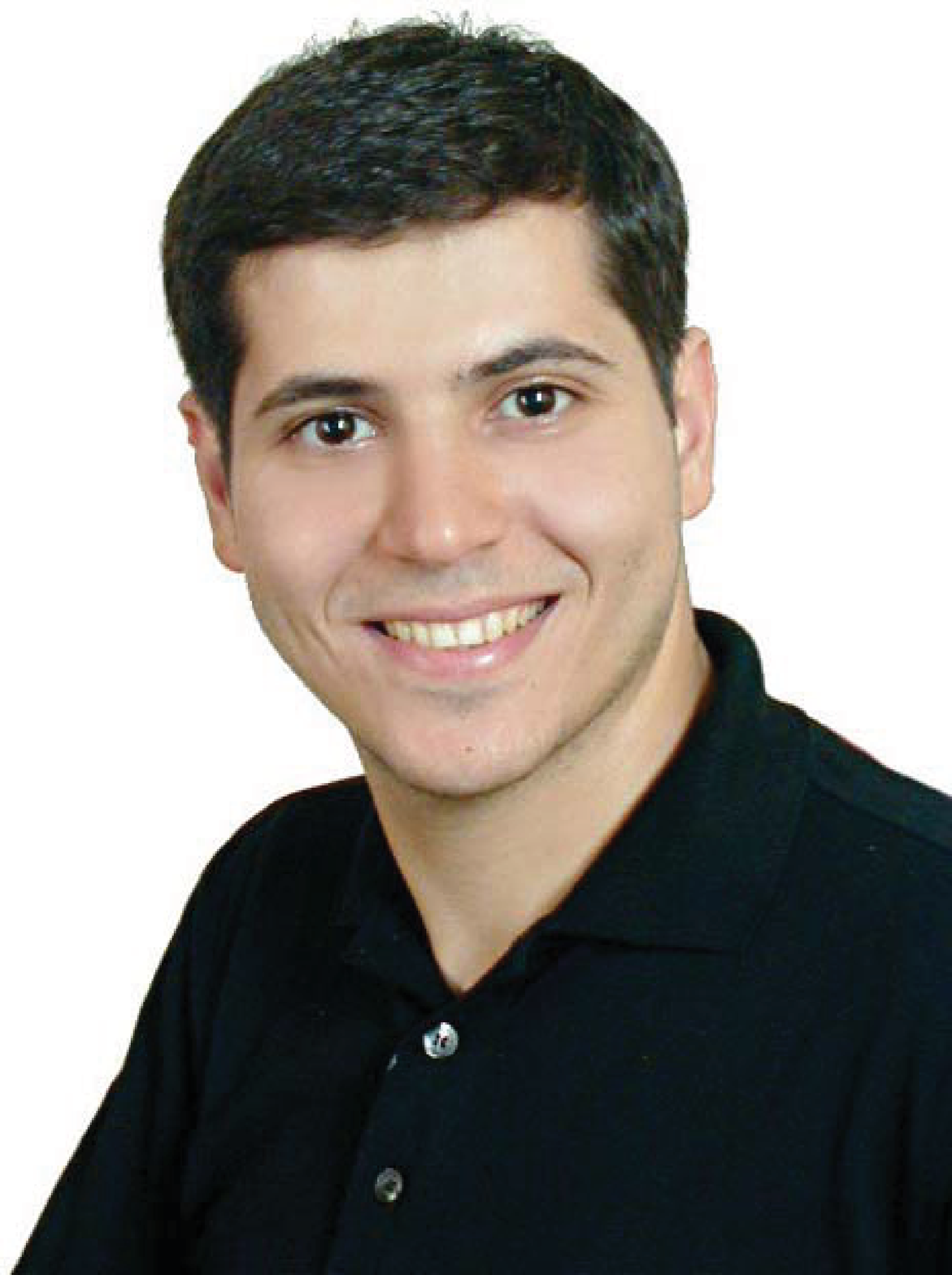}}]
{Zaher Dawy} 
 received the B.E. degree in Computer and Communications Engineering from the American University of Beirut (AUB) in 1998. He received his M.E. and Dr.-Ing. degrees in Communications Engineering from Munich University of Technology (TUM) in 2000 and 2004, respectively. He joined the Department of Electrical and Computer Engineering at AUB in September 2004 where he is currently an Associate Professor. Dr. Dawy is the recipient of the AUB 2008 teaching excellence award, best graduate award from TUM in 2000, youth and knowledge Siemens scholarship for distinguished students in 1999, and distinguished graduate medal of excellence from Harriri foundation in 1998. He is a senior member of the IEEE, Chair of the IEEE Communications Society Lebanon Chapter, and a member of the Lebanese Order of Engineers. His research interests are in the general areas of computational biology, information theory, and wireless communications with focus on genomic coding theory, gene network modeling, distributed and cooperative communications, cellular technologies, radio network planning and optimization, and multimedia transmission over communication networks.
\end{IEEEbiography}

\begin{IEEEbiography}[{\includegraphics[width=1in,height=1.25in,clip,keepaspectratio]{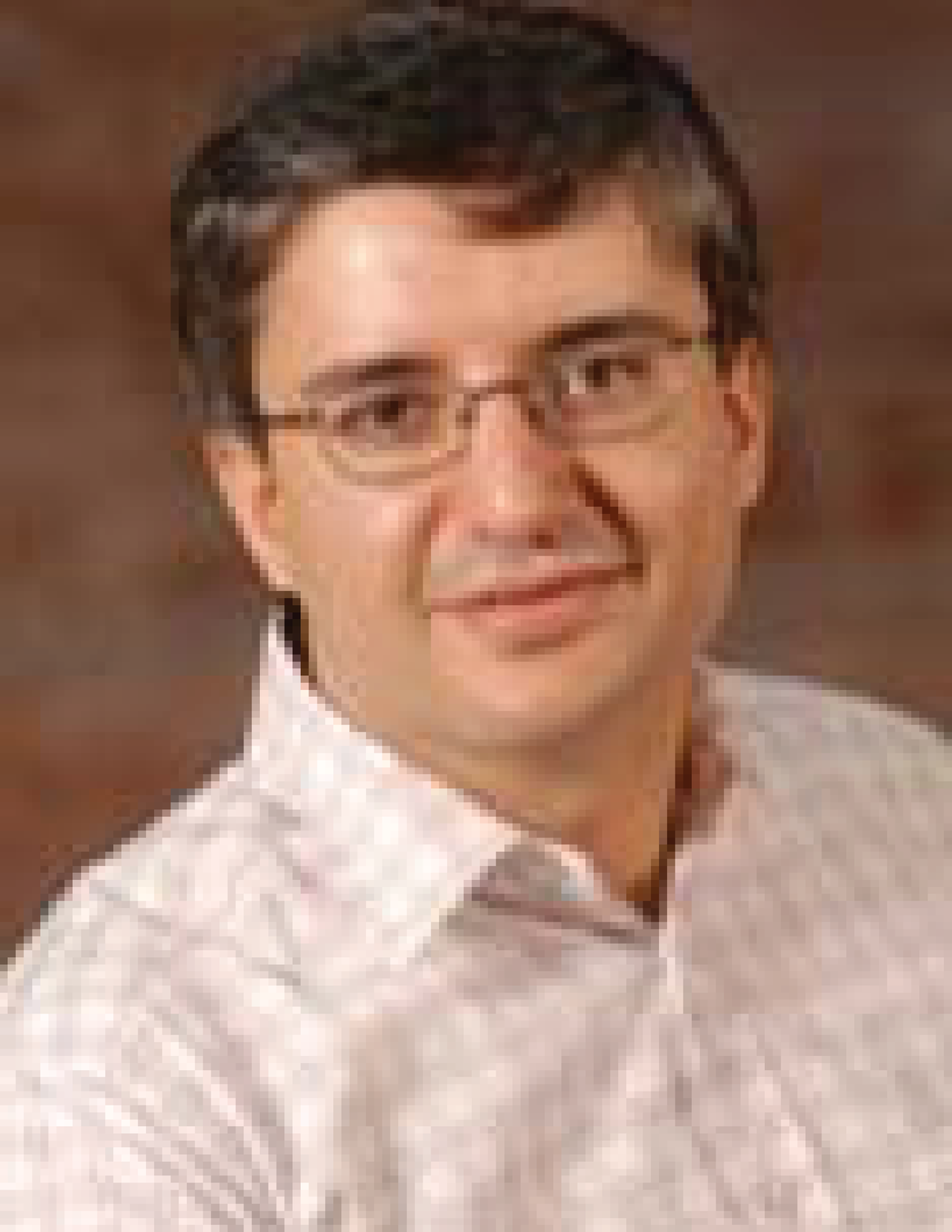}}]{Mohamed-Slim Alouini}(S'94, M'98, SM'03, F'09) was born in Tunis, Tunisia. He received the Ph.D. degree in electrical engineering from the California Institute of Technology (Caltech), Pasadena, CA, USA, in 1998. He was with the department of Electrical and Computer Engineering of the University of Minnesota, Minneapolis, MN, USA, then with the Electrical and Computer Engineering Program at the Texas A\&M University at Qatar, Education City, Doha, Qatar. Since June 2009, he has been a Professor of Electrical Engineering at King Abdullah University of Science and Technology (KAUST), Makkah Province, Saudi Arabia, where his current research interests include the modeling, design, optimization, and performance analysis of wireless communication systems.
\end{IEEEbiography}

\end{document}